\documentclass[a4paper,10pt]{article}

\usepackage[latin1]{inputenc}
\usepackage[francais,english]{babel}
\usepackage{graphicx,url}
\usepackage{latexsym,amssymb,amsmath,amsthm}
\usepackage[T1]{fontenc}

\theoremstyle{definition}

\theoremstyle{remark}

\numberwithin{equation}{section}

\begin{document}

\title{A case of mathematical eponymy: the Vandermonde determinant}
\author{Bernard Ycart\\
Laboratoire Jean Kuntzmann\\ 
Universit\'e Joseph Fourier and CNRS UMR 5224\\ 
51 rue des Math\'ematiques 38041 Grenoble cedex 9, France \\
\texttt{Bernard.Ycart@imag.fr}}
\maketitle
\begin{abstract}
We study the historical process that led to the worldwide adoption,
throughout mathematical research papers and textbooks, of
the denomination ``Vandermonde determinant''. The
mathematical object can be related to two passages in
Vandermonde's writings, of which one inspired Cauchy's definition of
determinants. Influential citations of Cauchy and Jacobi may have
initiated the naming process. It started during the second half of the
19\textsuperscript{th} century as a teaching practice in France. 
The spread in textbooks and
research journals began during the first half of
20\textsuperscript{th} century, and only reached full 
acceptance after the 1960's. The naming
process is still ongoing, in the sense that the volume of publications
using the denomination grows significantly faster than the
overall volume of the field.
\end{abstract}
\selectlanguage{francais}
\begin{abstract}
Nous \'etudions le processus historique qui a conduit \`a l'adoption
dans le monde entier de la d\'enomination \og d\'eterminant de
Vandermonde\fg{}. L'objet math\'ematique peut \^etre reli\'e \`a deux
passages dans les \'ecrits de Vandermonde, dont l'un a
inspir\'e Cauchy pour sa d\'efinition des d\'eterminants. Les
citations de Cauchy et Jacobi ont pu d\'eclencher le processus de
d\'enomination. Celui-ci a d\'emarr\'e au cours de la seconde moiti\'e
du \textsc{xix}\ieme{} si\`ecle comme une pratique p\'edagogique. Cette
pratique a pr\'ec\'ed\'e, plut\^ot que suivi, la diffusion dans les
livres et les articles de recherche, qui a commenc\'e pendant la
premi\`ere moiti\'e du \textsc{xx}\ieme{} si\`ecle, et n'a atteint un
r\'eel consensus qu'apr\`es les ann\'ees 1960. Le processus de
d\'enomination est encore en cours, au sens o\`u l'usage du nom
cro\^{\i}t significativement plus vite que le volume global de
publications du domaine.  
\end{abstract}
\selectlanguage{english}

{\small

}
\section{Introduction}
The Vandermonde determinant has become a standard example of
Stigler's law of eponymy: ``No scientific discovery is named after its
original discoverer'' (see \cite[p.~277]{Stigler1999}). 
The source? An authority: Henri Lebesgue
(1875--1941). On October 20 1937, 
he gave a conference at Utrecht University, entitled 
``L'{\oe}uvre math\'ematique de Vandermonde''. The
text of that conference was published in 1939, reproduced in 
1956, and again in a  1958 monography to which we shall 
refer: \cite{Lebesgue1958}. In order
to enhance Vandermonde's main achievement on the
resolution of algebraic equations \cite{Vandermonde1770}, Lebesgue
downplays his three other memoirs
\cite[p.~21]{Lebesgue1958}\footnote{all translations are B.Y.'s}:
\begin{quotation}
Thus the Vandermonde determinant is not due to Vandermonde; his theory
of determinants is not very original, his notation of factorials is
unimportant; his study of situation geometry is somewhat childish,
what is left? 
\end{quotation}
Actually, the memoir on combinatorics \cite{Vandermonde1772a} 
contains more than just a
notation for factorials: the identity
$$
\binom{n}{k}=\sum_{i=0}^k \binom{m}{i}\binom{n-m}{k-i}
$$
is still referred to as ``Vandermonde's theorem'' in probability and
combinatorics textbooks
(e.g. p.~315 of \cite{Santos2011}). Though ``childish'', 
the memoir on situation geometry \cite{Vandermonde1772b}
made him regarded as a precursor of knot theory (see 
\cite{Przytycki1992}).
\vskip 2mm
The life of Alexandre Th\'eophile Vandermonde (1735--1796), 
his engagement during the French revolution, 
his interests in music, mechanics, and political economy, and
his short mathematical carrier, have all been amply documented: see
\cite{Lebesgue1958}, \cite{Hecht1971}, \cite{Gillispie1976},
\cite{Faccarello1993}, and \cite{Sullivan1997}. We
shall not attempt a new biography nor a mathematical assessment of
Vandermonde's contribution. Neither shall we review here the early 
history of determinants. T. Muir's
\emph{Theory of determinants in their historical order of development}
is the indispensable basis, and we shall often refer to the first
two volumes: \cite{Muir1906,Muir1911}. Our focus here is exclusively on the
Vandermonde determinant, and more precisely on how
that particular object came to be known under that name. We call
Vandermonde Determinant, and denote by VD hereafter, the following
determinant, depending on $n$ variables $a_1,\ldots,a_n$.
\begin{equation*}
\left|
\begin{array}{ccccc}
1&a_1&a_1^2&\ldots&a_1^{n}\\
1&a_2&a_2^2&\ldots&a_2^{n}\\
\vdots&\vdots&\ddots&&\vdots\\
\vdots&\vdots&&\ddots&\vdots\\
1&a_n&a_n^2&\ldots&a_n^{n}
\end{array}
\right|
\end{equation*}
The VD has different mathematically equivalent interpretations, as a
product of differences or an alternating polynomial, that
will be developed in section \ref{cauchy}.
\vskip 2mm
Lebesgue makes the following 
assertion \cite[p.~21]{Lebesgue1958}:
\begin{quotation}
What could have been personal, is the Vandermonde determinant?
  Yet it is not there, nor anywhere else in Vandermonde's work.
\end{quotation}
Why then was Vandermonde's name given to that determinant? Lebesgue
has a conjecture.
\begin{quotation}
Vandermonde considers linear equations of which the unknowns are
denoted by $\xi 1, \xi 2, \xi3, \ldots$, and the coefficient of
$\xi i$ in the $k$\textsuperscript{th} equation by
$\stackrel{\scriptstyle{k}}{i}$. The resolution of such a system, 
e.g. of
$$
\begin{array}{lcl}
\stackrel{\scriptstyle{1}}{1}\xi 1+
\stackrel{\scriptstyle{1}}{2}\xi 2+
\stackrel{\scriptstyle{1}}{3}\xi 3+
\stackrel{\scriptstyle{1}}{4}&=&0,\\[1.5ex]
\stackrel{\scriptstyle{2}}{1}\xi 1+
\stackrel{\scriptstyle{2}}{2}\xi 2+
\stackrel{\scriptstyle{2}}{3}\xi 3+
\stackrel{\scriptstyle{2}}{4}&=&0,\\[1.5ex]
\stackrel{\scriptstyle{3}}{1}\xi 1+
\stackrel{\scriptstyle{3}}{2}\xi 2+
\stackrel{\scriptstyle{3}}{3}\xi 3+
\stackrel{\scriptstyle{3}}{4}&=&0,
\end{array}
$$   
will give determinants such as
$$
\left|
\begin{array}{ccc}
\stackrel{\scriptstyle{1}}{1}&
\stackrel{\scriptstyle{1}}{2}&
\stackrel{\scriptstyle{1}}{3}\\[1.5ex]
\stackrel{\scriptstyle{2}}{1}&
\stackrel{\scriptstyle{2}}{2}&
\stackrel{\scriptstyle{2}}{3}\\[1.5ex]
\stackrel{\scriptstyle{3}}{1}&
\stackrel{\scriptstyle{3}}{2}&
\stackrel{\scriptstyle{3}}{3}
\end{array}
\right|~;
$$
Forgetting for a while the convention of notations that has been made,
interpret the upper indices as exponents, you get a Vandermonde
determinant. And perhaps, it is this mix-up 
that saves Vandermonde's name from a more complete oblivion.
\end{quotation}
As we shall see, no trace of such a mix-up can be found in the
literature. Quite on the contrary, the mutation of
exponents into indices in a VD
is the very foundation of Cauchy's theory of determinants
\cite{Cauchy1812b}. Vandermonde himself \cite{Vandermonde1771}
had made the observation that changing one of the indices of a general 
determinant into an exponent led to an alternating function. That
remark did not escape either Cauchy
nor Jacobi; this may have
been the most solid argument in favor of the naming. On the other
hand, it does not quite make the VD
a counterexample to Stigler's law: 
linear systems with Vandermonde matrices had been written and
solved long before Vandermonde, by Isaac Newton (1642--1727)
and Abraham de Moivre (1667--1754).
\vskip 2mm
Nevertheless, our purpose here is not to decide 
whether it is right or wrong to name that determinant after 
Vandermonde (the reader will be given 
enough elements to form his/her own opinion). Neither is it to enter
the debate on mathematical eponymy (see \cite{Henwood1980,Smith1980}).
The naming of the VD is taken as a fact; and 
the \emph{history of that fact}, we believe, is of independent
interest. A mere attribution (citation: ``a determinant
introduced by Vandermonde'') must be distinguished
from an actual naming (eponymy: ``a Vandermonde
determinant''). The respective roles of citation (as a moral norm) and
eponymy (as a reward) in the
sociology of science have long been separated, following the pioneering
studies of R.K. Merton (e.g. \cite{Merton1968}). 
We refer to \cite{Small2004} for 
different theories of citation in science, and to \cite{Beaver1976}
for a historical perspective on eponymy. Eponymy has evolved together 
with successive sociological practices of science. In mathematics,
it became a widespread habit essentially during the
19\textsuperscript{th} century. 
Relatively few studies have been devoted to
mathematical eponymy; among them Stigler's articles (see
\cite{Stigler1999} and references therein) stand out.
The naming of a mathematical notion is in many cases
a long term process that extends over several generations of
mathematicians, and can be traced through historical
accounts, textbooks, and research publications. By \emph{naming process}
we mean the penetration of the name as a function of time,
``penetration'' being taken in the marketing sense: the proportion of
mathematicians knowing or using the name, measured as a proportion of
texts where it can be found.
\vskip 2mm
Lebesgue addressed his 1937 audience as follows.
\begin{quotation}
[\ldots] the name of Vandermonde would be ignored by the vast majority
of mathematicians if it had not been attributed to the determinant that
you know well, and which is not his!
\end{quotation}
The sentence seems to imply that the denomination ``Vandermonde
determinant'' was familiar to any mathematics student or professor in
1937. We believe that the naming process started as a teaching
practise  during the
second half of the 19\textsuperscript{th} century in France. 
Initially, it was more like a rumor than an identified
decision grounded on historical facts; actually, many mathematicians
clearly resisted it. As \cite[p.~283]{Stigler1999}
expresses it:
\begin{quotation}
[\ldots] resistance to eponymic recognition of close associates may
in fact be the norm of scientific behavior, one which serves the role
of protecting the practice from degenerating to a regional or
factional basis, with the consequent fall in the reward's incentive power.
\end{quotation}
This raises the question of the differential penetration of the naming
according to the countries, and the possible influence of
nationalisms, which we did not try to assess; it may be the case that
in 1937 the denomination was more familiar to Lebesgue than to his Dutch
audience.
The naming process of the VD
slowly gained momentum during the first half of the
20\textsuperscript{th}, but the denomination became
universally used by mathematicians only after the 1960's. It may be
considered that the naming process is still ongoing, in the 
sense that its
growth rate remains higher than that of the field. 
\vskip 2mm
To support
our assertions, we have examined a selection of influential textbooks,
conducted a systematic search through available databases, and
statistically studied numerical output data from MathSciNet. 
The first pedagogical publication we could find using the
denomination, appeared in 1886; the first textbook in 1897; the first
research paper in 1914.  
We have made a systematic query for the expressions ``Vandermonde
determinant'' and ``Vandermonde matrix'', on the MathSciNet database.
The occurrences start in 1929 and remain quite sporadic
until 1960. After 1960, the numbers of occurrences grow
exponentially. We have compared the growth rate with that of the (much
larger) number of occurrences of ``determinant'' or ``matrix''.
A statistical test has shown that the growth rate for  ``Vandermonde
determinant'' or matrix is significantly higher than the global
rate of increase for determinant or matrix. With all necessary
precautions on the use of quantitative methods (see
\cite{Goldstein1999}), our conclusion is that the
naming process, far from being an immediate recognition of
Vandermonde's achievements, 
is a rather recent, and still developing phenomenon. 
It appears to be posterior, and related, to 
the spread of matrix theory (see \cite{Brechenmacher2010}).
\vskip 2mm
The paper is organized as follows. Section \ref{alternants} 
gives a historical sketch of the mathematical
objects under consideration (difference-products and alternating functions).
Vandermonde's notations will be briefly examined in
\ref{notation}, 
then Cauchy's definition of determinants, based on difference-products, will be
exposed in \ref{cauchy}. 
In \ref{newton}, Newton's and de Moivre's anteriority
on the Vandermonde matrix through the divided differences method will
be reviewed. In
\ref{V123}, Vandermonde's actual contributions will  be
discussed. 
Section \ref{process} is devoted to the naming process, that will be
examined from three different points of view.
Historical accounts will be described in \ref{historical}, focusing on
the credits explicitly given to Vandermonde.
The appearance of the naming in textbooks is described in
\ref{textbooks}. The quantification of the
naming process in research papers 
is treated in \ref{databases}.
\section{Difference-products and alternating functions}
\label{alternants}
\subsection{Vandermonde's notation}
\label{notation}

Before describing the mathematical objects under study, 
we shall briefly
comment on Vandermonde's notations, of which Lebesgue thought they could
have induced a mix-up between indices and exponents. 
Here is Vandermonde's definition of determinants
\cite[p.~517]{Vandermonde1771}:
\begin{quotation}
I suppose that one represents by
$\stackrel{\scriptstyle{1}}{1}$, 
$\stackrel{\scriptstyle{2}}{1}$,
$\stackrel{\scriptstyle{3}}{1}$, \&c.
$\stackrel{\scriptstyle{1}}{2}$, 
$\stackrel{\scriptstyle{2}}{2}$,
$\stackrel{\scriptstyle{3}}{2}$, \&c.
$\stackrel{\scriptstyle{1}}{3}$, 
$\stackrel{\scriptstyle{2}}{3}$,
$\stackrel{\scriptstyle{3}}{3}$, \&c.
as many different general quantities,
of which any one be 
$\stackrel{\scriptstyle{\alpha}}{a}$,
another one be 
$\stackrel{\scriptstyle{\beta}}{b}$, \&c.
\& that the product of both be ordinarily denoted by
$\stackrel{\scriptstyle{\alpha}}{a}\cdot\stackrel{\scriptstyle{\beta}}{b}$.
Of the two ordinal numbers $\alpha$ \& $a$, the first one, for
instance, will designate from which equation the coefficient
$\stackrel{\scriptstyle{\alpha}}{a}$ is taken, and the second one will
designate the rank that the coefficient has in the equation, as will
be seen hereafter.

I suppose moreover the following system of abbreviations, and that
it be set  
{\small
$$
\begin{array}{l}
\begin{array}{c|c}\alpha&\beta\\\hline a&b\end{array}
\;=\;
\stackrel{\scriptstyle{\alpha}}{a}\cdot\stackrel{\scriptstyle{\beta}}{b}
-
\stackrel{\scriptstyle{\alpha}}{b}\cdot\stackrel{\scriptstyle{\beta}}{a}
\\[3ex]
\begin{array}{c|c|c}\alpha&\beta&\gamma\\\hline a&b&c\end{array}
\;=\;
\stackrel{\scriptstyle{\alpha}}{a}
\cdot\;
\begin{array}{c|c}\beta&\gamma\\\hline b&c\end{array}
\;+
\stackrel{\scriptstyle{\alpha}}{b}
\cdot\;
\begin{array}{c|c}\beta&\gamma\\\hline c&a\end{array}
\;+
\stackrel{\scriptstyle{\alpha}}{c}
\cdot\;
\begin{array}{c|c}\beta&\gamma\\\hline a&b\end{array}
%\\[3ex]
%\begin{array}{c|c|c|c}\alpha&\beta&\gamma&\delta\\\hline a&b&c&d\end{array}
%\;=\;
%\stackrel{\scriptstyle{\alpha}}{a}
%\cdot\;
%\begin{array}{c|c|c}\beta&\gamma&\delta\\\hline b&c&d\end{array}
%\;-
%\stackrel{\scriptstyle{\alpha}}{b}
%\cdot\;
%\begin{array}{c|c|c}\beta&\gamma&\delta\\\hline c&d&a\end{array}
%\;+
%\stackrel{\scriptstyle{\alpha}}{c}
%\cdot\;
%\begin{array}{c|c|c}\beta&\gamma&\delta\\\hline d&a&b\end{array}
%\;-
%\stackrel{\scriptstyle{\alpha}}{d}
%\cdot\;
%\begin{array}{c|c|c}\beta&\gamma&\delta\\\hline a&b&c\end{array}
\end{array}
$$
}
[\ldots]
\end{quotation}
Vandermonde's notations probably looked much less
strange in the 19\textsuperscript{th} century than
they do nowadays. 
Referring to them, T. Muir said \cite[p.~24]{Muir1906}:
\begin{quotation}
[\ldots] we observe first that Vandermonde proposes for coefficients
 a positional notation essentially the same as that of Leibnitz[sic], 
writing
 $\stackrel{\scriptstyle{1}}{2}$ where Leibnitz wrote $12$ or $1_2$.
\end{quotation}
Indeed, Vandermonde's notations were quite similar to some of the many
systems tried by Leibniz (see \cite{Knobloch2001}). During the first
half of the 19\textsuperscript{th} century, different ways of denoting
the coefficients of an array or a linear system coexisted
(see \cite{Muir1906}): ${^ij}$,
$i_j$, $(i,j)$, ${^ia_j}$, $a_i^{(j)}$\ldots~
In the first treatise ever published on determinants, W. Spottiswoode 
used $(i,j)$ \cite{Spottiswoode1851}.  
 C.L. Dodgson (Lewis Carroll)  was the only one who ever denoted
 coefficients by
$i\!\big\rmoustache\! j$ \cite{Dodgson1867}.
G. Dostor's classical treatise \cite{Dostor1877}
proposed two notations, ``juxtaposed'' and
``superposed'' indices. Suarez and Gasc\'o 
describe and use 6 different notations
\cite{Suarez1882}.
The modern notation $a_{i,j}$ was already present in Cauchy's 
memoir \cite[p.~113]{Cauchy1812b}. But Cauchy himself mostly preferred
multiple letter notations such as $a_i,b_i,c_i,\ldots,e_i,f_i$
(e.g. p.~121).  
\subsection{Cauchy's definition}
\label{cauchy}
Cauchy's two founding memoirs \cite{Cauchy1812a,Cauchy1812b} were read
to the Institute on November 30 1812, but were only published in 1815.
After a thorough analysis of both, T. Muir concludes with a very lively
description of the respective roles of Vandermonde
and Cauchy \cite[p.~131]{Muir1906}.
\begin{quotation}
If one bears this in mind, and recalls the fact, temporarily set
aside, that Cauchy, instead of being a compiler, presented the subject
from a perfectly new point of view, added many results previously
unthought of, and opened up a whole avenue of fresh investigation, one
cannot but assign to him the place of honour among all the workers
from 1693 to 1812. It is, no doubt, impossible to call him, as some
have done, the formal founder of the theory. This honour is certainly
due to Vandermonde, who, however, erected on the foundation
comparatively little of a superstructure. Those who followed
Vandermonde contributed, knowingly or unknowingly, only a stone or
two, larger or smaller, to the building. Cauchy relaid the foundation,
rebuilt the whole, and initiated new enlargements; the result being an
edifice which the architects of to-day may still admire and find
worthy of study.
\end{quotation}
What was that  ``perfectly new point of view''?
Previously, B\'ezout, Laplace, and Vandermonde 
had all defined determinants by induction 
using, explicitly or not, what is now known as Laplace's formula: the
development of 
a determinant along one of its lines or columns.  
Cauchy's definition  \cite[p.~113]{Cauchy1812b} is radically different:
\begin{quotation}
Let $a_1, a_2,\ldots,a_n$ be several different quantities in number
equal to $n$. It has been shown above, that by multiplying the product
of these quantities, or
$$
a_1a_2a_3\ldots a_n
$$
by the product of their respective differences, or else by
$$
(a_2-a_1)(a_3-a_1)\ldots(a_n-a_1)(a_3-a_2)\ldots(a_n-a_2)\ldots(a_n-a_{n-1})
$$
one obtained as a result the alternating symmetric function
$$
S(\pm a_1a_2^2\ldots a_n^n)
$$
which, as a consequence, happens to be always equal to the product
$$
a_1a_2\ldots a_n
(a_2-a_1)\ldots(a_n-a_1)(a_3-a_2)\ldots(a_n-a_2)
\ldots(a_n-a_{n-1})\;.
$$
Let us suppose now that one develops this later product and that,
in each term of the development, one replaces the exponent of each
letter by a second index equal to that exponent: by writing, for
instance, $a_{r,i}$ instead of $a_r^i$ and $a_{i,r}$ instead of
$a_i^r$, one will obtain as a result a new alternating symmetric
function which, instead of being represented by
$$
S(\pm a_1^1a_2^2\ldots a_n^n)\;,
$$
will be represented by
$$
S(\pm a_{1,1}a_{2,2}\ldots a_{n,n})\;,
$$
the sign $S$ being relative to the first indices of each letter.
Such is the most general form of the functions that I shall designate
in what follows under the name of \emph{determinants}.
\end{quotation}
In order to understand Cauchy's reasoning, one must keep in mind 
that his main focus was on functions of $n$ variables: 
\cite{Cauchy1812b} came as a sequel to \cite{Cauchy1812a} where he
discussed functions 
of $n$ variables that take less than $n!$ different values when the
variables are permuted. 
He called ``symmetric alternating functions'' (fonctions sym\'etriques
altern\'ees) 
those functions taking only two opposite values (they will be referred
to as ``alternating functions''). 
Among them, the polynomials in $n$ variables are multiples of the
``product of differences'', later called difference-product 
(see \cite{Muir1906}). The difference-product develops into a sum of
monomials with  
alternating signs. Those signs depend on the permutation of the
variables and their exponents,  
and the ``rule of signs'' had been described by Cauchy before defining
determinants. 
(On the discovery by Leibniz in 1683 of the rule of signs, see
\cite{Knobloch2001}). 
\vskip 2mm
Different expressions in $n$ variables $a_1,a_2,\ldots,a_n$, may both be
mathematically equivalent, and have different interpretations. We shall
distinguish between:
\begin{itemize}
\item difference-product: 
$\displaystyle{\prod_{1\leqslant i<j\leqslant n} (a_j-a_i)}$, 
\item alternating polynomial: 
$\displaystyle{\sum_{\sigma\in\mathcal{S}_n}
  (-1)^{\varepsilon(\sigma)}\prod_{i=1}^n a_i^{\sigma(i)-1}}$, 
\item Vandermonde determinant:
  $\mbox{det}(a_i^j)_{1\leqslant i\leqslant n, 0\leqslant j \leqslant n-1}$.
\end{itemize}
They are written in modern notations: $\mathcal{S}_n$ is the 
group of permutations of $\{a_1,\ldots,a_n\}$ onto itself and
$\varepsilon(\sigma)$ denotes the 
signature of the permutation $\sigma$. Needless to say, the group of
permutations  and the signature as a homomorphism are anachronistic. Cauchy
had recognized in the development  
of the difference-product, the same rule of signs as that of a general
determinant. Hence his idea of using
$$
\prod_{i=1}^n a_i \prod_{1\leqslant i<j\leqslant n}^n (a_j-a_i)
=
\sum_{\sigma\in\mathcal{S}_n} (-1)^{\varepsilon(\sigma)}\prod_{i=1}^n a_i^{\sigma(i)}
$$
as a general definition, after mutating the exponent of
each variable into a second index. 
\vskip 2mm\noindent
As pointed out by \cite[p.~247]{Muir1906}, the year 1841 marked a new
spurt for determinant 
theory, fueled by the publication in Crelle's journal of Jacobi's
monograph, divided into 3 papers. 
There Jacobi rebuilds the whole theory, taking Cauchy's approach upside down.
Here is Muir's account \cite[p.~254]{Muir1906}.
\begin{quotation}
At the outset, there is a reversal of former orders of things; Cramer's rule 
of signs for a permutation and Cauchy's rule
being led up by a series of propositions instead of one of them being
made a convention or definition. This 
implies, of course, that a new definition of a signed permutation is
adopted, and that conversely this definition 
must have appeared as a deduced theorem in any exposition having
either of this rules as its starting point. 
\end{quotation}
In other words, when Cauchy's started from the difference-product,
then defined a general determinant 
by mutating exponents into indices, Jacobi first defined positive and
negative permutations, then  
defined the determinant as a polynomial, with coefficients 
$\pm 1$ according to the sign
of the permutation.
Eventually, Jacobi's definition prevailed upon Cauchy's, which was 
forgotten.  
Cauchy undoubdtedly saw both pedagogical and mathematical advantages
to his approach. When he writes his famous 
``Cours d'Analyse'' in 1821, he follows exactly the same path as in
his 1812 memoir. He recommends  
the difference-product as a general method for solving linear systems
of equations, and applies it immediately 
to the Lagrange interpolation problem (pp.~71, 72, 426, 429 of
\cite{Cauchy1821}). The third of Jacobi's memoirs 
in Crelle's Journal \cite{Jacobi1841} deals with alternating
functions. Cauchy responds 
with \cite{Cauchy1841} in which he treats quotients of alternating functions
by difference-products.  
In particular, he calculates the determinant 
$\mathrm{det}\left(\frac{1}{a_i+b_j}\right)_{1\leqslant i,j\leqslant n}$ (formula
(10) p.~154 of \cite{Cauchy1841}) in a quite simple way. 
(Interestingly enough, the denomination 
``Cauchy determinant'' for that example seems to be rarely used
outside France, whereas the particular case 
$a_i=i, b_j=j-1$ is universally known as ``Hilbert matrix'').   
\vskip 2mm
One year before 1841, the difference-product approach had been
rediscovered by James Joseph Sylvester (1814-1897) 
\cite{Sylvester1840}, who (without any reference to Cauchy) 
called ``zeta-ic multiplication'' Cauchy's operation of mutating
exponents into indices in a polynomial. 
Muir's comment \cite[p.~235]{Muir1906} is somewhat ironic.
\begin{quotation}
This early paper, one cannot but observe, has all the characteristics
afterwards so familiar to readers of Sylvester's writings, -- fervid
imagination, 
vigorous originality, bold exuberance of diction, hasty if not contemptuous
disregard of historical research, the outstripping of demonstration by
enunciation, and an infective enthousiasm as to the vistas opened by
his own work. 
\end{quotation} 
\subsection{Newton, de Moivre, and the interpolation problem}
\label{newton}
The difference-product could hardly be considered an original notion
in Cauchy's time.  
Apart from being a very natural way of combining $n$ variables, it
appears in the Lagrange 
interpolation  problem. This other interesting case of mathematical eponymy 
is connected to ours, as we shall now see. For a history of interpolation, see
\cite{Fraser1919}, and section 10.4 of \cite{Chabert1999}.
If $(x_1,y_1),\ldots,(x_n,y_n)$ are the Cartesian coordinates of the
points to be interpolated 
and $P=a_0+a_1x+\cdots+a_{n-1}x^{n-1}$ the unknown polynomial, then
its coefficients
$a_0,\ldots a_{n-1}$ satisfy the following linear system.
\begin{equation*}
\tag{LIS}
\left\{\begin{array}{ccc}
a_0+a_1x_1+\cdots+a_{n-1}x_1^{n-1}&=&y_1\\
a_0+a_1x_2+\cdots+a_{n-1}x_2^{n-1}&=&y_2\\
\vdots&\vdots&\vdots \\
a_0+a_1x_n+\cdots+a_{n-1}x_n^{n-1}&=&y_2 
\end{array}\right.
\end{equation*}
Assuming the $x_i$'s are all different,
the solution is the Lagrange interpolation polynomial:
\begin{equation*}
\tag{LIP}
P(X)=\sum_{i=1}^n y_i \prod_{j\neq i}\frac{X-x_j}{x_i-x_j}\;.
\end{equation*}
It may seem fair that whoever first wrote the system of equations
(LIS) should get the credit 
for discovering the Vandermonde matrix and whoever wrote (LIP)
for computing its inverse (and implicitly the VD).  
The naming ``Lagrange interpolation''
comes from one of the lessons that Joseph Louis Lagrange (1736--1813)
gave at the \'Ecole Normale in Paris in 1795 \cite{Lagrange1795}. 
There, Lagrange did not pretend to expose
his own research:
\begin{quotation}
Newton is the first one who has posed that problem. Here is the
solution he gives. [\ldots]
\end{quotation}
Indeed, in 
the \emph{Principia Mathematica}, 
Isaac Newton (1642--1727) had described a method to determine ``a
curved line of parabolic type 
which passes through any number of points'' 
\cite[pp.~695--696]{Newton1687}: 
what is now known as
Newton's divided differences 
method. In the \emph{Principia}, Newton did not explicitly write (LIS). 
However, in a famous letter to Oldenburg
dated October 24 1676, he mentions a manuscript, \emph{Methodus
differentialis}, that appeared in print only 
after the \emph{Principia}, in 1711. There, the system (LIS) is explicitly 
written (see
p.~10 of \cite{Fraser1919}, where 
the \emph{Methodus Differentialis} is reproduced and translated), but
the explicit solution (LIP) is not given. One may think that writing down (LIP)
would have seemed useless and even misleading to Newton: he must have
been aware
that his method was both faster and numerically
more stable than the direct application of (LIP).
The first one to explicitly write (LIP)  is Newton's friend
Abraham de Moivre (1667--1754), in 1730 (on de Moivre's relationship
with Newton, see \cite{Bellhouse2007}). 
Instead of interpolation, de Moivre's motivation was to calculate the
coefficients in a  
linear combination of geometric series, when that linear combination
is supposed equal  
to another series. The coefficients turn out to be the solution of a
system equivalent to (LIS). 
In theorem \textsc{iv} pp.~33--35 of the \emph{Miscellanea analytica} 
\cite{Moivre1730},
de Moivre explicitly writes a 
general system with  power coefficients, and gives its solution, thus
being the first one to write 
the inverse of a Vandermonde matrix. Actually, de Moivre had already
published particular cases of that result 
in the first edition of 
his \emph{Doctrine of chances} 
\cite[p.~132]{Moivre1718}. There he said:
\begin{quotation}
And if a general theorem were desired, it might easily be formed from
the inspection 
of the foregoing.

These theorems are very useful for summing up readily those series
which express the probability of 
the plays being ended in a given number of games.
\end{quotation}
Indeed, de Moivre's motivation came from probability problems arising
from dice games: 
the theorem is used for the solution of problem \textsc{iv}, p.~77 of
\emph{Miscellanea analytica}, and in later editions 
(1738 and 1756) of the 
\emph{Doctrine of chances}. 
De Moivre gives full credit to Newton both for the interpolation
problem and the divided differences method.  
The following extract of his preface to the \emph{Doctrine of chances}
\cite[p.~\textsc{x}]{Moivre1718} is worth quoting: its last sentence
has a particular resounding with our subject.
\begin{quotation}
There are other sorts of series, which tho' not properly infinite, yet
are called series, from the  
regularity of the terms whereof they are composed; those terms
following one another with a certain 
uniformity, which is always to be defined. Of this nature is the
Theorem given by Sir \emph{Isaac 
Newton}, in the fifth \emph{Lemma} of the third Book of his
\emph{Principles}, for drawing a curve through 
any given number of points: of which the demonstration, as well as
other things belonging to the same subject, 
may be deduced from the first Proposition of his \emph{Methodus
  Differentialis}, printed with some other of his  
tracts, by the care of my intimate friend, and very skilful
mathematician, Mr. W. Jones. The abovementionned 
theorem being very useful in summing up any number of terms whose last
differences are equal (such as the 
numbers called triangular, pyramidal, \&c. the squares, the cubes, or
other powers of numbers in arithmetic  
progression) I have shewn in many places of this book how it might be
applicable to these cases. I hope it will 
not be taken amiss that I have ascribed the invention of it to its
proper author, tho' 'tis possible 
some persons may have found something like it by their own sagacity.
\end{quotation}
De Moivre's anteriority on the difference-product has been
pointed out on several occasions, in particular by 
\cite{Tee1993}; but of course, de Moivre does not express
difference-products as determinants.  Actually, the difference-product, 
and the explicit expression 
of the inverse matrix have been rediscovered many times, until late in
the 20\textsuperscript{th} century: see \emph{e.g.}
\cite{Klinger1967}. 
\subsection{Vandermonde's writings}
\label{V123}
We shall now examine what in Vandermonde's work can be connected to
the VD.  
About his ``Memoir on elimination'', Vandermonde says
(\cite{Vandermonde1771}, footnote p.~516):
\begin{quotation}
This memoir was read to the Academy for the first time on the
20\textsuperscript{th} 
of January 1771. It contained different things that I have suppressed
here because they have been published 
since by other Geometers.
\end{quotation}
These ``other Geometers'' certainly include Laplace, whose memoir
though posterior, was
published in the same  
volume as Vandermonde's. Guessing what exactly did Vandermonde 
suppress cannot but remain
conjectural.  
\vskip 2mm
Just like Cauchy in 1812, Vandermonde wrote about determinants as a
byproduct of symmetric functions; 
his memoir on elimination is a sequel to the
memoir on the resolution of equations. 
The publications dates, 1774 and 1776, are misleading:
\cite{Vandermonde1770} was read to the academy  
``sometime in November 1770'', \emph{i.e.} only two months before
\cite{Vandermonde1771}.  
Vandermonde undoubtedly had the first memoir in mind when he wrote the
second, and both should be examined 
as a whole. Here are two quotations, that we have numbered for later
reference. 
\begin{itemize}
\item[[V1\!\!]] \cite[p.~369]{Vandermonde1770}:
\begin{quotation}
And yet, $(a^2b+b^2c+c^2a-a^2c-b^2a-c^2b)$, which equals
$(a-b)(a-c)(b-c)$, squares as
$$
\begin{array}{c}
a^4b^2+a^4c^2+b^4c^2+c^4a^2+c^4b^2\\
-2(a ^4bc+b^4ac+c^4ab)-2(a^3b^3+a^3c^3+b^3c^3)\\
+2(a^3b^2c+a^3c^2b+b^3a^2c+b^3c^2a+c^3a^2b+c^3b^2a)\!-\!6a^2b^2c^2.
\end{array}
$$
\end{quotation}
\item[[V2\!\!]] \cite[p.~522]{Vandermonde1771}:
\begin{quotation}
Those acquainted with the abbreviated symbols that I have named
\emph{partial types of combination}, in my \emph{Memoir on the
resolution of equations}, will recognize here the formation
of the \emph{partial type} depending on the second degree, for any
number of letters; they will easily see that, by taking our
$\alpha, \beta, \gamma, \delta$, \&c. for instance, as exponents,
all terms with equal signs in the development of one of our
abbreviations, will also  
be the development of the \emph{partial type} depending on the second
degree, \& formed 
with an equal number of letters.
\end{quotation}
\end{itemize}
Actually, the difference-product of four variables appears in the
following passage \cite[p.~386]{Vandermonde1770}:
\begin{quotation}
The first of these cubes is
$$
\begin{array}{l}
(A^3B^3)-\frac{3}{2}(A^3B^2C)+6(A^3BCD)+6(A^2B^2C^2)-3(A^2B^2CD)\\
+\frac{3}{2}(a-b)(a-c)(a-d)(b-c)(b-d)(c-d)\sqrt{-3}\;;
\end{array}
$$
[\ldots] as the
square of the  
product of differences between the roots is a function of
\emph{types}, [\ldots]
\end{quotation}
However, the development is not explicitly written, 
and we have not found that sentence ever referred to.
\vskip 2mm
In \cite{Vandermonde1770}, Vandermonde details the resolution of
second and third degree equations 
(hence [V1]), then states his general method, and illustrates it by
the fourth degree equation. The rest of the paper is
devoted to a discussion on the symmetric functions of the 
roots. Admittedly,
the difference-product of three variables appears in [V1], and 
its development is given; but this does not establish 
that Vandermonde saw it as a determinant.
[V2] certainly proves that he knew  determinants were related to his 
``partial types depending on the second degree'' (\emph{i.e.}
alternating functions), 
through changing indices into exponents. He probably knew exactly to which
``partial type'' did the VD correspond, at least in
dimension 3, and probably 
in dimension 4. There is no evidence he actually wrote a VD
as a particular determinant, nor that he wrote difference-products 
of more than four variables. The impressive tables displayed
on the three pages after p.~374 of \cite{Vandermonde1770} show that he certainly
had the capacity for much more difficult formal calculations. But they
also prove that he did not have a general expression for symmetric nor
alternating functions. The long footnote of pp.~374--375 seems to
imply that he was  
on his way towards greater generality.
\begin{quotation}
[\ldots] 
By considering this formula as a multivariate finite difference
equation, in which 
the difference of each variable is equal to unity, I can integrate \&
satisfy the conditions, 
by a particular procedure of which I propose to render an account in one
of the future volumes. 
\end{quotation}
It is not very surprising that, by manipulating symmetric functions of
3 or 4 variables, Vandermonde 
had been led to write difference-products.  Whether or not he viewed
them as determinants may not be  
the most important. More interesting is the relation that he had seen
in [V2]. He undoubtedly knew that  
by making an exponent of the second index in a determinant, 
an alternating function
was obtained. But conversely,   
had he realized that \emph{any} determinant could be obtained from a
difference-product by the reverse operation? 
[V2] comes in \cite{Vandermonde1771}, immediately after his 4 pages 
``proof'' of the alternating property, before which he had announced:
\begin{quotation}
Instead of generally proving these two equations [the alternating property],
which would demand an awkward rather than difficult calculation, I
shall content myself 
with developing the simplest examples; this will suffice to grasp the
spirit of the 
proof.
\end{quotation}
The alternating property of the difference-product is trivial; 
and with Cauchy's definition, proving that a determinant changes sign
when exchanging two co\-lumns 
becomes obvious. We do not think that 
Vandermonde would have written his four pages of
``simplest examples'' 
had he anticipated Cauchy's definition. Lebesgue appreciation 
 on Vandermonde's
contribution to the resolution of equations 
 might still have some truth in
it when applied to  
Vandermonde's determinants \cite[p.~38]{Lebesgue1958}:
\begin{quotation}
Vandermonde never came back on his algebraic researches because at first
he felt only imperfectly their importance, and if he did not
understand it better afterwards, it is precisely 
because he had not reflected deeply on them; [\ldots]
\end{quotation}

\section{The naming process}
\label{process}
\subsection{Historical accounts}
\label{historical}
We have searched historical notes in textbooks or research papers, for
connections being made between Vandermonde and the VD. 
Many accounts have been given of
Vandermonde's contribution to the resolution of equations: see 
\cite{Neumann2007} or \cite{Stedall2011} for recent
references. Among the most famous, \cite{Nielsen1929} and 
\cite{Waerden1985} (as many others) do not mention the VD. 
Similarly Vandermonde's founding role is
acknowledged in most historical accounts of determinant theory, but
there again, his relation to the VD is seldom mentioned: 
throughout history, there seems to have been some embarrassment on
the subject.

Muir's masterly treatise is quite significant,
and it may have had some later influence on the naming.
As many other authors, Muir calls ``difference-product'' the VD and
``alternants'' those determinants stemming from alternating functions or
generalizing the VD; he has been quite an active contributor
of the field in the last decades of the 19\textsuperscript{th}
century. In each volume, he devotes a chapter to alternants. Here are
the first lines of that chapter in Volume 1
\cite[p.~306]{Muir1906}:
\begin{quotation}
The first traces of the special functions now known as
\emph{alternating functions} 
are said by Cauchy to be discernible in certain work of Vandermonde's;
and if we view the functions as 
originating in the study of the number of values which a function can
assume through permutation 
of its variables, such an early date may in a certain sense be
justifiable. To all intents and purposes, however, the 
theory is a creation of Cauchy's, and it is almost absolutely certain
that its connection with determinants was never thought 
of until his time.
\end{quotation}
In volume 2, Muir feels obliged to set some records straight
\cite[p.~154]{Muir1911}, p.~154:
\begin{quotation}
Further, as exagerated statements regarding Vandermonde's contribution to the
subject have been widely accepted, it seems desirable to point out the
exact foundation on which such statements
rest. In a paper read in November 1770 Vandermonde says (p.~369), ``Or
$a^2b+b^2c+c^2a-a^2c-b^2a-c^2b$, qui \'egale $(a-b)(a-c)(b-c)$ a pour carr\'e
$a^4b^2+\ldots$'' This is the whole matter.
\end{quotation}
As we have seen, there are essentially two ways to connect
Vandermonde's writings to the VD:
\begin{itemize}
\item[[V1\!\!]]: Vandermonde has written the difference-product of three
  variables and its development, hence a particular case of the VD.
\item[[V2\!\!]]: Vandermonde has anticipated Cauchy's definition by
  remarking that changing one of the indices into an exponent gives an
  alternating function.
\end{itemize}
Clearly, Muir is on the [V1] side, as all historians have been
since. It was not quite so in the
19\textsuperscript{th} century. As Muir points out,
Cauchy had studied Vandermonde's two memoirs on the resolution of
equations and on elimination, and quotes them. In 
\cite[p.~110]{Cauchy1812b}, [V1] is explicitly cited:
\begin{quotation}
Thus, supposing for instance $n=3$, it will be found 
\begin{eqnarray*}
S^2(\pm a_2a_3^2)&=&
a_2a_3^2+a_3a_1^2+a_1a_2^2-a_3a_2^2-a_2a_1^2-a_1a_3^2\\
&=&(a_2-a_1)(a_3-a_1)(a_3-a_2)\;.
\end{eqnarray*}
This last equation has been given by Vandermonde in his memoir on the
resolution of equations.
\end{quotation}
Cauchy does not explicitly acknowledge that [V2] 
inspired his definition of determinants
from difference-products, but the following quotation
clearly alludes to [V2] \cite[p.~70]{Cauchy1812a}.
\begin{quotation}
The smallest divisor of this product is equal to $2$ and it is
easy to make sure, that, in any order, it is possible to form
functions having only two different values. Vandermonde has given ways
to compose functions of that kind. In general, to form with quantities
$$
a_1,a_2,\ldots,a_n
$$
an order $n$ function with index $2$, it will suffice to consider the
positive or the negative part of the product
$$
(a_1-a_2)(a_1-a_3)\cdots(a_1-a_n)(a_2-a_3)\cdots(a_2-a_n)\cdots(a_{n-1}-a_n)
$$
whose factors are the differences of the quantities 
$a_1,a_2,\ldots,a_n$ taken two by two.
\end{quotation}
\vskip 2mm
We could find in the literature only 4 other citations of [V2]. The
earliest comes in the very first words of \cite{Jacobi1841}; 
admittedly, it is worth many others.
\begin{quotation}
The famous Vandermonde once elegantly observed that the proposed
determinant
$$
 \sum \pm a_0^{(0)}a_1^{(1)}a_2^{(2)}\ldots a_n^{(n)}\;,
$$
if indices are changed into exponents, comes from the product formed
from the differences of all elements $a_0,a_1,\ldots,a_n$
$$
\begin{array}{r}
P=(a_1-a_0)(a_2-a_0)(a_3-a_0)\cdots(a_n-a_0)\\
(a_2-a_1)(a_3-a_1)\cdots(a_n-a_1)\\
(a_3-a_2)\cdots(a_n-a_2)\\
\cdots\cdots\\
(a_n-a_{n-1})
\end{array}
$$
\end{quotation}
The next citation that we are aware of, appears in 
\cite{Terquem1846}.
\begin{quotation}
A very ingenious observation of the same geometer [Vandermonde], 
about
indices considered as exponents, has given birth to Mr. Cauchy's
beautiful theory of \emph{alternating functions} and to his 
proof of Cramer's 
formulae.
\end{quotation}
Our third citation comes from the preface of
Spottiswoode's treatise. There he 
comments \cite{Cauchy1812b} as follows
\cite[p.~vi]{Spottiswoode1851}: 
\begin{quotation}
The second part of this paper refers immediately to determinants,
and contains a large number of very general theorems. Amongst
them is noticed a property of a class of functions closely connected
with determinants, first given, so far as I am aware, by Vandermonde;
if in the development of the expression
$$
  a_1a_2\cdots a_n
  (a_2-a_1)\cdots(a_n-a_1)(a_3-a_2)\cdots(a_n-a_2)\cdots(a_n-a_{n-1})
$$
the indices be replaced by a second series of suffixes, the 
result will be the determinant
$$
S(\pm a_{1,1}a_{2,2}\ldots a_{n,n})\;.
$$
\end{quotation}
The last citation appears in \cite{Prouhet1856} who, before writing
the difference-product of $n$ variables ``according to a theorem due
to Vandermonde'' gives [V2] as a reference.
\vskip 2mm 
It is likely that, since Cauchy's definition never prevailed
and soon fell into oblivion, so went with it Vandermonde's ``elegant
observation''. From then on, [V1] has been the commonly
accepted source for the naming. The position usually adopted is clearly
expressed by an anonymous contributor to the ``Nouvelles annales de
Math\'ematiques'' \cite[p.~181]{Anonymous1860}.
\begin{quotation}
Vandermonde [\ldots] decomposes into factors a polynomial that can be
considered as a 3\textsuperscript{rd} order determinant: but nothing
indicates that he had the general theorem in mind, not even that he
had considered that polynomial as a determinant.
\end{quotation}
The same view has been expressed many times, from R. Baltzer
\cite[p.~50]{Baltzer1857} to J. Stedall \cite[p.~190]{Stedall2011}, 
through S. G\"unther \cite[p.~66]{Gunther1875}
and G. Kowalewski \cite[p.~315]{Kowalewski1942}; it appears in the
Encyclopedia of Mathematics \cite[p.~363]{Remeslennikov1993}.
Only two of the early authors were less carefull in their attribution:
F. Brioschi speaks of an ``important relation due to
Vandermonde'' \cite[p.~75]{Brioschi1854}, and 
G.A. Gohierre de Longchamps devotes a section to ``Vandermonde's
theorem''
\cite[p.~82]{Gohierre1883}. 
\vskip 2mm
Since the publication of Lebesgue's conference \cite{Lebesgue1958}, 
his mix-up conjecture has been cited by several 
authors: see \emph{e.g.}
\cite[p.~18]{Edwards1984}, \cite[p.~197]{Blyth2002};
it even appears in 
Gillispie's Dictionary of Scientific Biography, \cite[p.~571]{Gillispie1976}.
It has probably fostered the widely accepted idea that the
attribution of the VD to Vandermonde is a misnomer. 
J. Dieudonn\'e states it quite clearly \cite[p.~59]{Dieudonne1978}.
\begin{quotation}
This naming, due to Cauchy, is not historically justified, since
Vandermonde never explicitly introduced such a determinant.
\end{quotation}
Yet, Dieudonn\'e was aware of Cauchy's use of the exchange between
exponents and indices, that he presents as an ``elegant
trick''\ldots{}
\subsection{Textbooks}
\label{textbooks}
We have made a selection of 24 treatises and
textbooks having appeared in the
19\textsuperscript{th} and 20\textsuperscript{th} centuries, 
partially or completely devoted to determinants, and where the VD
appears as a mathematical object, if only as a simple example or
exercise. All of them have had several
editions or translations, which we regard as a criterion of
(relatively) large diffusion. Our selection is
arbitrary, and we have examined only a very small
sample of the full textbook production of these times.
We have not systematically
searched outside the area of linear algebra, though we are aware that
early occurrences of the naming can be found in other fields. For
instance,
in one of the earliest and most influential treatises on numerical
analysis, when the authors expose Newton's divided
difference method, they 
write the interpolation system, its determinant, and add 
\cite[p.~23]{Whittacker1924}:
\begin{quotation}
Now a difference-product may be expressed as a determinant of the kind
known as Vandermonde's[\ldots]
\end{quotation}
As another example, P\'olya and Szeg\H{o}'s famous textbook contains a
``generalized Vandermonde determinant'' \cite[p.~43]{Polya1945}.
Nevertheless, we
consider our sample as representative, in the statistical sense: our
conclusion being that the denomination remains sporadic until 1950,
we believe it would be confirmed on a broader corpus. 
Table \ref{tab:books} gives the references, the publication
country (including translations), and the name given to the VD for
each book in our sample.
\begin{table}[!ht]
\begin{tabular}{|l|l|c|l|}
\hline
Reference&Countries&Page&Naming\\\hline
\cite{Brioschi1854}&Italy&75&none\\
\cite{Baltzer1857}&Germany, France&50&none\\
\cite{Salmon1859}&Great-Britain&13&none\\
\cite{Bertrand1859}&France, Italy&333&none\\
\cite{Trudi1862}&Italy&31&none\\
\cite{Gunther1875}&Germany&66&difference product\\
\cite{Dostor1877}&France&142&none\\
\cite{Scott1880}&Great-Britain&115&difference product\\
\cite{Mansion1880}&Belgium&27&none\\
\cite{Suarez1882}&Spain&360&none\\
\cite{Gohierre1883}&France&82&none\\
\cite{Hanus1886}&USA&187&difference product\\
\cite{Chrystal1886}&Great-Britain&53&none\\
\cite{Pascal1897}&Italy, Germany&166&Vandermonde\\
\cite{Kronecker1903}&Germany&304&none\\
\cite{Hawkes1905}&USA&218&none\\
\cite{Weld1906}&USA&169&alternant\\
\cite{Wedderburn1934}&USA&26&none\\
\cite{Barnard1936}&Great-Britain, USA&126&none\\
\cite{Aitken1939}&USA&42&alternant\\
\cite{Kowalewski1942}&Germany&315&none\\
\cite{Gantmacher1953}&Russia, USA&99&Vandermonde\\
\cite{Bourbaki1970}&France, USA&532&Vandermonde\\
\cite{Lang1970}&USA&155&Vandermonde\\
\hline
\end{tabular}
\caption{Textbooks including the VD, and whether ot not it is given a
  name.}
\label{tab:books}
\end{table}
\vskip 2mm
Before the second half of the 20\textsuperscript{th} century, the
denomination  ``Vandermonde determinant'' 
can hardly be found in textbooks.
Among the early treatises on determinants, \cite[p.~75]{Brioschi1854}
mentions ``an important relation due to Vandermonde'', and
\cite{Gohierre1883} devotes a section to ``Vandermonde's
theorem''. These attributions 
may have had some influence on the naming practice, 
but they are not \emph{actual namings} of the VD as a
mathematical object.
Ernesto Pascal (1865--1940) seems to be the first one to actually name the
VD in a textbook. His hesitations are very revealing. The running head of 
\cite[p.~166]{Pascal1897}
is indeed ``Vandermonde determinant''. 
But the title of the section is ``Vandermonde or
Cauchy determinant''. Pascal cites \cite{Jacobi1841} and mentions:
\begin{quotation}
It is usually called also Cauchy determinant, this last author
having considered it in general, 
whereas Vandermonde studied it in a
particular case. 
\end{quotation}
Many authors, although quite aware of Vandermonde's contributions,
remain very cautious regarding the naming.
Siegmund G\"unther (1848--1923) 
devotes the first chapter of his treatise to a careful historical
exposition, where Vandermonde's role is thoroughly analyzed. 
Yet later on, the VD
is named ``Differenzenproduct'' and attributed to Vandermonde 
for $n=3$ and
to Cauchy for the general case \cite[p.~66]{Gunther1875}.
Leopold Kronecker (1823--1891) cannot be suspected of downplaying
Vandermonde's achievements (see \cite{Lebesgue1958}). 
However, when he writes his 
``Lessons on the theory of
determinants'', he attributes the VD to Cauchy 
\cite[p.~304]{Kronecker1903} and does not name it. In his ``Lessons
on number theory'', the VD is named ``Differenzenprodukt''
\cite[p.~396]{Kronecker1901}.
Joseph Bertrand (1822--1900) has known Cauchy, and he 
is among the rare authors to follow
Cauchy's definition of determinants. His ``Trait\'e \'el\'ementaire 
d'alg\`ebre'' had several editions since 1851. The
determinants appear in the 1859 Italian edition
 \cite[p.~333]{Bertrand1859} but no name is given to the VD.
\subsection{Research papers}
\label{databases}
\begin{figure}[!ht]
\centerline{
\includegraphics[width=10cm]{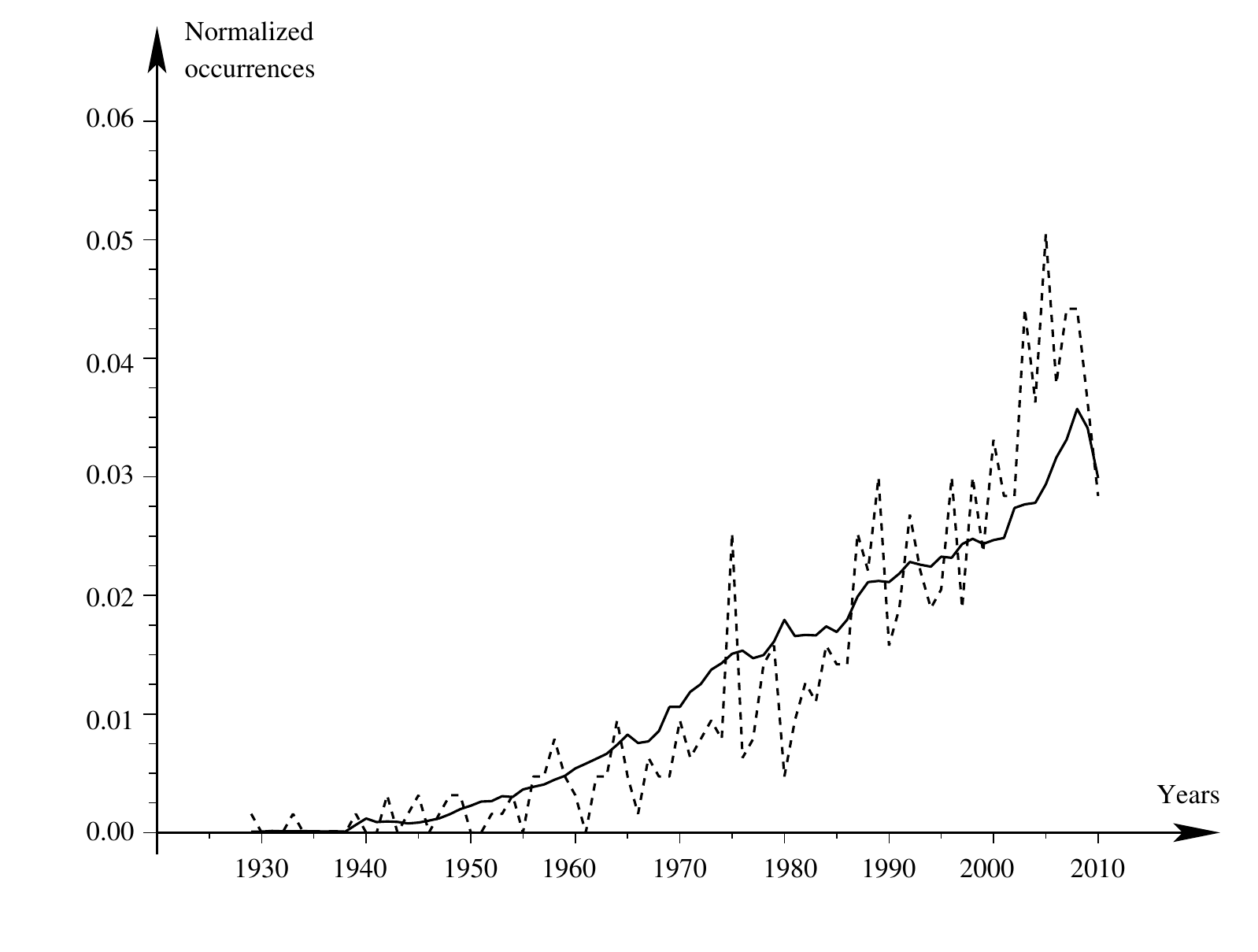}
} 
\caption{Occurrences of ``Vandermonde determinant'' or ``Vandermonde
  matrix'' (dashed) compared to ``determinant'' or ``matrix'' (solid) in the
  MathSciNet database. For each curve, the data per year have been 
divided by their sum.}
\label{fig:MSN1}
\end{figure}

\begin{figure}[!ht]
\centerline{
\includegraphics[width=10cm]{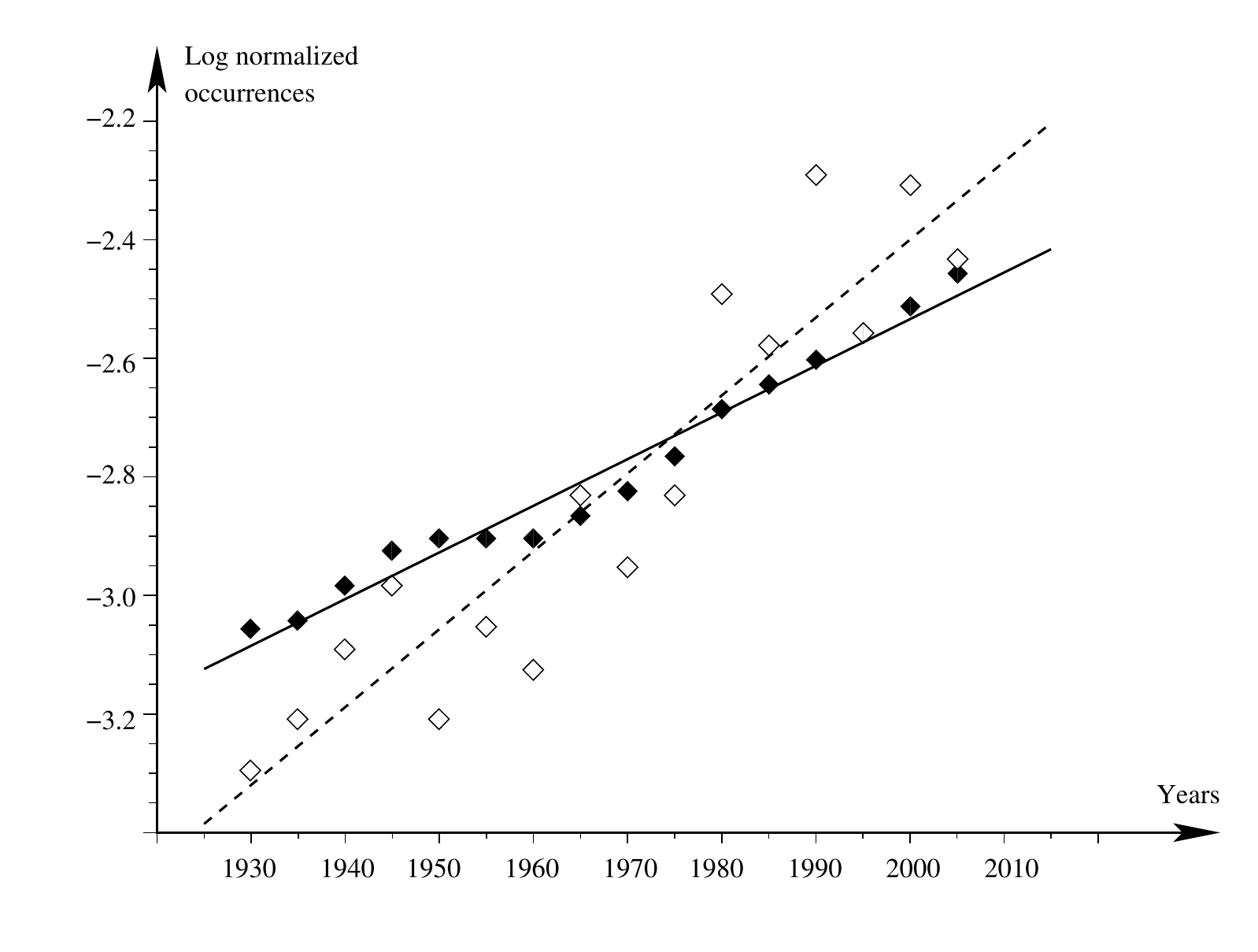}
} 
\caption{Linear regressions for the logarithms of
occurrences of ``Vandermonde determinant'' or ``Vandermonde
  matrix'' (dashed line, empty diamonds) and ``determinant'' or ``matrix'' 
  (solid line and diamonds) in the
  MathSciNet database. The data are binned by 5-year periods over
  the 80 years 1929-2008.}
\label{fig:MSN2}
\end{figure}

In order to evaluate the penetration of the expression ``Vandermonde
determinant'' in the mathematical literature, we have
searched through several data\-ba\-ses: Gallica, Google Books,
G\"ottinger Digitalisierungszentrum, Internet Ar\-chive, Jstor,
Mathematical Reviews (or ``MathSciNet''), Numdam, and
Zentralblatt Math\footnote{
\url{http://gallica.bnf.fr},
\url{http://books.google.com},
\url{http://gdz.sub.uni-goettingen.de},
\url{http://www.archive.org},
\url{http://www.jstor.org},
\url{http://www.ams.org/mathscinet},
\url{http://www.numdam.org/},
\url{http://www.zentralblatt-math.org/zbmath}.}.
The earliest traces of the attribution 
that we could find in articles are:
\begin{enumerate}
\item
\cite[p.~87]{Prouhet1856}: ``According to a theorem due to Vandermonde'' 
\item
\cite[p.~181]{Anonymous1860}: ``This theorem, ordinarily attributed
to Vandermonde,[\ldots]''
\item
\cite[p.~517]{Neuberg1866}: ``This last determinant, by virtue of
the theorem known as \emph{Vandermonde's},[\ldots]''
\end{enumerate}
We cannot be sure that earlier appearances do not exist elsewhere. 
However we find it significant that
the earliest references were found in pedagogy rather than research
journals. They come from professors at the undergraduate level, 
sharing their solutions to particular problems. 
In quotations 2 and 3, some hesitation
can be felt in the expressions ``ordinarily attributed to'' or ``known
as''. As we have already seen, \cite{Prouhet1856} cites [V2] to support
the attribution, whereas \cite{Anonymous1860} clearly resists
it; both implicitly admit that the attribution to
Vandermonde is already a usual practice.
After 1886, maybe under the influence of
\cite{Gohierre1883}, the attributions become more assertive.
The first two actual namings seem to be:
\begin{enumerate}
\item
\cite[p.~164]{Marchand1886}: ``The numerator is a Vandermonde
determinant.''
\item
\cite{Weill1888}~: ''On a form of Vandermonde determinant'' (title of
the paper).
\end{enumerate}
The first occurrence of the naming in a research journal was found through
Jstor: \cite{Bennett1914}. This indicates that the denomination was
already in use both among researchers and outside
France,
before World War One. 
\vskip 2mm
For our quantitative study, we chose to focus on MathSciNet, 
that seemed to give more easily
interpretable results. As an example of the difficulties encountered
with other bases, 
Zentralblatt has references to which the
keyword ``Vandermonde determinant'' is associated, whereas it does
not appear in the article: an example is \cite{Jonquieres1895} whose
denomination for the VD is ``d\'eterminant potentiel'';
these false detections were difficult to sort. However we believe
that searching in another database would give similar results (compare 
Figure \ref{fig:MSN1} below with those of Annex 1.2 in
\cite{Brechenmacher2010}).
We are aware of the
limits to our quantitative approach. The MathSciNet database does not
contain all published articles; moreover, we could not check
each reference to make sure it was relevant. Nevertheless, we
consider that MathSciNet is a representative sample, in the
statistical sense, of the
total mathematical production: we believe that our estimation of
exponential growth rates would not be significantly (again in
the statistical sense) modified if computed on another database.
\vskip 2mm
We first
searched for the other historical denominations, ``alternant'',
``diffe\-rence-product'' and 
``power determinant''. No publication
could be found for ``power determinant'', which seems to have disappeared 
(maybe for ambiguity reasons). Similarly, only two non ambiguous
occurrences were found for  ``difference-product''.
The name ``alternant'' is also ambiguous: 
it appears in ``alternant code'' and ``alternant
group''. After disambiguation, here are the occurrences per decade.
$$
\begin{array}{|c|cccccccc|}
\hline
\mbox{dates}&<1940&40\mbox{'s}&50\mbox{'s}&60\mbox{'s}&70\mbox{'s}&
80\mbox{'s}&90\mbox{'s}&>2000\\
\hline
\mbox{alternant occurrences}&10&7&8&12&4&9&5&4\\
\hline
\end{array}
$$   
The occurrence of ``alternant'' (as a determinant) 
did not
completely disappear, but it has remained sporadic, and has not
increased with the total mathematical production.
Let us now turn to the Vandermonde denomination. It can be
found under different forms.
\begin{itemize}
\item Vandermonde determinant or matrix,
\item Vandermonde's determinant or matrix,
\item Vandermondian.
\end{itemize}
The second one has 16 occurrences before 2011, the third one only 7.
The first occurrence of ``Vandermondian'' was found in
\cite{Farrel1959}; however, the term seems 
to be more current in the physical
literature than in the mathematical one: see \cite{Vein1999}, section
4.1 p.~51. It may be the case that the use of the Vandermonde
determinant in the modelling of the quantum Hall effect 
(see \cite{Scharf1994}) boosted its
popularity among physicists. This would match the effect that quantum
mechanics had on the development of matrix theory, as described by
\cite{Brechenmacher2010}.  
\vskip 2mm
The query ``Vandermonde determinant'' includes
``Vandermonde's determinant'' (and determinants); applied with the
option ``Anywhere'', it returns 273
occurrences. The query ``Vandermonde matrix'' (including plural)
returns 363 occurrences. Our query was the disjunction of these
two, and it returned 623 occurrences 
(less than the sum of the previous
two because ``determinant'' and ``matrix'' together are found in 13
references). The first occurrence appears in 1929. 
We have made the same query for each year from 1929 to 2010. The
corresponding numbers will be referred to as ``Vandermonde data''. 
They remain quite sporadic
during the first half of the 20\textsuperscript{th} century 
(0, 1, 2,
or 3 occurrences per year before 1958); then they gradually  
increase. 
Of course that increase was expected, since 
the total mathematical production grows
exponentially: the increase in the output of 
any given query should be considered only relatively to the
increase of the
total production in the field. 
For the same years (1929--2010), we
have made the query ``determinant or matrix''. The corresponding series
will be referred to as ``global data''. The total number was 
202219. In order to compare both series, we have plotted on the same
graphic (Figure \ref{fig:MSN1}), the Vandermonde and the global
data, after dividing each by its sum. 
Of course the Vandermonde data are more irregular; however, 
both curves seem to grow exponentially, with a higher rate for
the Vandermonde data. 
\vskip 2mm
In order to provide a statistical justification to the previous assertions, 
our treatment was the following. Firstly, the last two years 
(2009 and 2010) were truncated:
they show a decrease that we do not consider as significant; it is
probably due to the delay in entering new publications in the
base. Then the data were binned over periods of 5 years (to account for
sporadicity at the beginning of the Vandermonde series). 
Saying that the data grow exponentially means that they can be
adjusted by a function 
of the type $y=\exp(ax+b)$ where $x$ is a year, $y$ a number of
publication, and $a$  
is the exponential growth rate. Equivalently, the logarithm of the
data can be adjusted 
by a linear function of the years: $ax+b$. The parameters $a$ and $b$
were estimated by a  
least-squares linear regression of the log-data over the years
(see \emph{e.g.} chap.~14 of \cite{Utts2004} as a general reference). Figure
\ref{fig:MSN2} 
displays the graphical results of the two linear regressions. 
Both regressions were found to be significant, with respective p-values of
$3.6~10^{-12}$ and $3.1~10^{-7}$. 
The exponential growth rate 
(\emph{i.e.} the slope of the regression line) was found to be    
$0.0079$ for the global data,
and $0.0131$ for the Vandermonde data. In other words, the global
number of publications is multiplied by $\mathrm{e}^a\simeq 1.0079$,
or else increases 
by $0.79\%$ per year on average, whereas the Vandermonde data increase
by $1.31\%$.  
To test whether the $0.52\%$ observed difference between 
growth rates was significant, we used another linear regression, that
time on the logarithm of the \emph{ratios}, \emph{i.e.} 
on the difference of the two previous sets. The new slope is of course
the difference of the two previous ones, and 
it was found to be significantly positive, with a p-value of $6.9~10^{-4}$.
\vskip 2mm
Having shown that the denomination  ``Vandermonde determinant or matrix'' 
has a higher growth rate than ``determinant or matrix'' alone, the
question of the interpretation arises. Comparing exponential growth
rates may be a way of measuring the scientific dynamism of a
research field. A field with a faster growth than the global
production could be considered as booming; on the contrary a field 
with a lower growth rate would be seen as slowing down; 
among two fields, the more dynamic would be
the one with a significantly higher growth rate.
Here, the problem is different. The hypothesis of a higher
dynamics of research on the VD compared to the rest of linear algebra
can be ruled out:
the VD has long been an undergraduate-level basic tool rather than a
subject of research of its own. There remains two possible
explanations.
\begin{enumerate}
\item
The fields of research using the Vandermonde determinant or matrix 
as a tool, are more fertile than those using other determinants or
matrices.
\item
Mathematicians using a Vandermonde determinant or matrix tend more and
more to give it its usual name. 
\end{enumerate} 
We could not find any evidence supporting the first hypothesis, and we
believe that the occurrence of the VD as an object is no more
frequent in today's mathematical research than it was some decades
ago. The only explanation we find plausible is that when
mathematicians encounter a VD, they tend more and more to use 
the standard denomination, which has become a universally accepted shortcut. 
\section{Conclusion}
\label{conclusion}
In our study of the historical process that led to the worldwide adoption,
throughout mathematical research papers and textbooks, of
the denomination ``Vandermonde determinant'', we have established the
following points. Although Vandermonde is not the first
discoverer of the object, although he never expressed it in full generality,
there still exist two connections between his writings and the VD:
he has written down and developed the difference-product of 3 variables,
and he has observed that changing indices into exponents in a general
determinant gave an alternating function. 
Even if Vandermonde's calculation of the 3
variables difference-product was the only one eventually retained 
by historians, his second observation about changing exponents into indices
probably inspired Cauchy's definition of determinants, and was quoted
by Jacobi. Both may have sparked off the naming process. It
 started during the second half of the 19\textsuperscript{th} century,
 essentially as a teaching practice. For quite a long time, 
textbook and research paper authors resisted the
naming, for which no sufficient justification existed in their view.
The naming process eventually gained momentum during the second half of the 
20\textsuperscript{th} century and from then on, its penetration of the
mathematical community has been increasing. This was
proved by a statistical treatment of numerical data from the
MathSciNet database, that consisted in 
comparing the exponential growth rates of the naming to that of the
global production.
\vskip 2mm
Thus we believe that we have brought answers to 
the questions where?, when?, and how? The most important 
question may be the one we did not address: why? The sociological
explanation of eponymy as a reward, may not be the only one. We believe
that the pedagogical function of eponymy, which has been overlooked until
now, should be taken into account. Here are some of the questions
that would deserve an investigation. As
the computation of the VD became a classical exercise or example, did
the pressure to name it increase? More generally, 
do students
prefer a mathematician's name rather than an impersonal one? Is a
theorem easier to memorize when given a person's name? 
Does a mathematician necessarily transmit as a researcher 
the denominations he has learned as a student? 
Many questions remain to be asked, but we do not think that they are proper
to mathematics, nor that can be
answered by mathematicians alone: maybe the time has
come for a collaboration between specialists of mathematics, pedagogy, and
onomastics (see \emph{e.g.} \cite{Nuessel2011})\ldots
\subsection*{Acknowledgements}
The author is indebted to the two anonymous referees for important
remarks and helpful hints.
\bibliographystyle{decsi}
%\bibliography{Vandermonde}

\begin{thebibliography}{A}

\bibitem[Aitken 1939]{Aitken1939}
A.C. Aitken,
\textit{Determinants and matrices},
Oliver \& Boyd, Edinburgh, 1939.

\bibitem[Baltzer 1857]{Baltzer1857}
R. Baltzer, \textit{Theorie und Anwendung der Determinanten},
Hirzel, Leipzig, 1857.

\bibitem[Barnard \& Child 1936]{Barnard1936}
S. Barnard \& J.M. Child, \textit{Advanced algebra},
Macmillan \& Co., New-York, 1936

\bibitem[Beaver 1976]{Beaver1976}
D.D. Beaver,
Reflections on the natural history of eponymy and scientific law,
\textit{Social Studies of Science}, 6, 1976, pp.~89--98.

\bibitem[Bellhouse \& Genest 2007]{Bellhouse2007}
D.R. Bellhouse \& C. Genest,
Maty's biography of Abraham de Moivre,
translated, annotated and augmented,
\textit{Statistical Science}, 22(1), 2007, pp.~109--136.

\bibitem[Bennett 1914]{Bennett1914}
A.A. Bennett, An algebraic treatment of the theorem of closure,
\textit{Ann. Math.}, 2\textsuperscript{nd} ser., 16(1/4), 1914,
pp.~97--118.

\bibitem[Bertrand 1859]{Bertrand1859}
G. Bertrand,
\textit{Trattato di algebra elementare},
Le Monnier, Firenze, 1859.

\bibitem[Blyth \& Robertson 2002]{Blyth2002}
T.S. Blyth \& E.F. Robertson,
\textit{Further linear algebra},
Springer, New York 2002.

\bibitem[Bourbaki 1970]{Bourbaki1970}
\textit{Elements of Mathematics: Algebra 
\textsc{i}},
Springer, New York, 1989.

\bibitem[Brechenmacher 2010]{Brechenmacher2010}
F. Brechenmacher,
Une histoire de l'universalit\'e 
des matrices math\'ematiques,
\textit{Revue de Synth\`ese}, 131(4), 2010, pp.~569-603.

\bibitem[Brioschi 1854]{Brioschi1854}
\textit{La teorica dei determinanti, e le sue pricipali applicazioni},
Bizzoni, Pavia, 1854.

\bibitem[Cauchy 1812a]{Cauchy1812a}
A.L. Cauchy,
M\'emoire sur le nombre des valeurs qu'une fonction peut acqu\'erir
lorsqu'on y permute de toutes les mani\`eres possibles les quantit\'es
qu'elle renferme,
\textit{Journal de l'\'Ecole Polytechnique}, 
\textsc{xvii}\textsuperscript{e}
Cahier, Tome \textsc{x}, 1815, pp.~64--90.

\bibitem[Cauchy 1812b]{Cauchy1812b}
A.L. Cauchy,
M\'emoire sur les fonctions qui ne peuvent obtenir
que deux valeurs \'egales et de signes contraires par suite des
transpositions op\'er\'ees entre les variables qu'elles renferment,
\textit{Journal de l'\'Ecole Polytechnique}, 
\textsc{xvii}\textsuperscript{e}
Cahier, Tome \textsc{x}, 1815, pp.~91--169.

\bibitem[Cauchy 1821]{Cauchy1821}
A.L. Cauchy,
\textit{Cours d'Analyse de l'\'Ecole Royale Polytechnique}
Debure, Paris, 1821.

\bibitem[Cauchy 1841]{Cauchy1841}
A.L. Cauchy,
M\'emoire sur les fonctions altern\'ees et sur les sommes
  altern\'ees,
in \textit{Exercices d'Analyse et de Physique Math\'ematique} 
tome \textsc{ii}, 
Bachelier, Paris 1841, pp.~151--159.

\bibitem[Chabert \& Barbin 1999]{Chabert1999}
J.L. Chabert \& \'E. Barbin,
\textit{A history of algorithms: from the pebble to the microchip},
Springer, New York, 1999.

\bibitem[Chrystal 1886]{Chrystal1886}
G. Chrystal, \textit{Algebra},
Black, London, 1886.

\bibitem[de Jonqui\`eres 1895]{Jonquieres1895}
E. de Jonqui\`eres,
Sur les d\'ependances mutuelles des d\'eterminants
potentiels,
\textit{CRAS Paris} 120, 1895, pp.~408--410.

\bibitem[de Moivre 1718]{Moivre1718}
A. de Moivre,
\textit{The doctrine of chances},
London, 1718.

\bibitem[de Moivre 1730]{Moivre1730}
A. de Moivre,
\textit{Miscellanea analytica de 
seriebus et quadraturis},
London, 1730.

%\bibitem[Dhombres 1998]{Dhombres1998}
%J. Dhombres,
%Une histoire de l'objectivit\'e scientifique 
%et le concept de post\'erit\'e,
%in R. Guesnerie \& H. Hartog, \textit{Des sciences et des techniques, un 
%d\'ebat}, EHESS, Paris 1998, pp.~127--148.

\bibitem[Dieudonn\'e 1978]{Dieudonne1978}
J. Dieudonn\'e,
\textit{Abr\'eg\'e d'histoire des math\'ematiques
1700--1900},
Tome \textsc{i}, Hermann, Paris, 1978.

\bibitem[Dodgson 1867]{Dodgson1867}
C.L. Dodgson,
\textit{An elementary treatise on determinants},
Mac Millan, London, 1867.

\bibitem[Dostor 1877]{Dostor1877}
\textit{\'El\'ements de la th\'eorie des d\'eterminants},
Delagrave, Paris, 1877.

\bibitem[Edwards 1984]{Edwards1984}
H.M. Edwards,
\textit{Galois theory},
Springer, New York, 1984.

\bibitem[Faccarello 1993]{Faccarello1993}
G. Faccarello,
Du conservatoire \`a l'\'Ecole Normale~:
quelques notes sur A.T. Vandermonde (1735--1796),
\textit{Cahiers d'Histoire du CNAM}, 2/3 1993, pp.~15--57.

\bibitem[Farrel 1950]{Farrel1959}
A.B. Farrel, A special Vandermondian determinant,
\textit{Amer. Math. Monthly}, 66, 1959, pp.~564--569.

\bibitem[Fraser 1919]{Fraser1919}
D.C. Fraser, Newton's interpolation formulas,
reprinted from \textit{The Journal of the Institute of Acturaries},
vol. li, pp.77--106 (Oct. 1918) and pp.~211--232
(April 1919). Layton, London, 1919.

\bibitem[Gantmacher 1953]{Gantmacher1953}
F.R. Gantmacher,
\textit{The theory of matrices},
Chelsea, New York, 1959.

\bibitem[Gillispie 1976]{Gillispie1976}
C.C. Gillispie,
\textit{Dictionary of scientific biography,
vol. \textsc{xiii}},
Scribner, New York, 1976.

\bibitem[Gohierre 1883]{Gohierre1883}
G.A. Gohierre de Longchamps
\textit{Cours de Math\'ematiques sp\'e\-cia\-les},
Delagrave, Paris, 1883.

\bibitem[Goldstein 1999]{Goldstein1999}
C. Goldstein,
Sur la question des m\'ethodes quantitatives en
histoire des math\'ematiques~: le cas de la th\'eorie des nombres en
France (1870--1914),
\textit{Acta historiae rerum naturalium nec non technicarum} 
3(28), 1999, pp.~187--214.

\bibitem[G\"unther 1875]{Gunther1875}
S. G\"unter,
\textit{Lehrbuch der Determinanten-theorie f\"ur Studirende},
Besold, Erlangen 1875.

\bibitem[Hanus 1886]{Hanus1886}
P.H. Hanus,
\textit{An elementary treatise on the theory of determinants},
Ginn \& Co., Boston, 1886.

\bibitem[Hawkes 1905]{Hawkes1905}
H.E. Hawkes, \textit{Advanced algebra},
Gin \& Co., Boston, 1905.

\bibitem[Henwood \& Rival 1980]{Henwood1980}
M.R. Henwood \& I. Rival,
Eponymy in Mathematical nomenclature:
what's in a name, and what should be?,
\textit{Math. Intelligencer}, 2(4), 1980, pp.~204--205.

\bibitem[Hecht 1971]{Hecht1971}
J. Hecht,
Un exemple de multidisciplinarit\'e~: Alexandre Vandermonde,
\textit{Population}, 26(4), 1971, pp.~641--676.

\bibitem[Jacobi 1841]{Jacobi1841}
C.G. Jacobi,
De fonctionibus alternatibus earumque divisione per productum
  et differentiis elementorum conflatum,
\textit{J. Reine Angew. Mathematik}, 
vol. \textsc{xxvii} 1841, pp.~360--371.

\bibitem[Klinger 1967]{Klinger1967}
The Vandermonde matrix,
\textit{Amer. Math. Monthly}, 74(5), 1967, pp.~571--574.

\bibitem[Kowalewski 1942]{Kowalewski1942}
G. Kowalewski,
\textit{Determinantentheorie},
de Gruyter, Berlin, 1942.

\bibitem[Knobloch 2001]{Knobloch2001}
E. Knobloch,
D\'eterminants et \'elimination chez Leibniz,
\textit{Revue d'Histoire des Sciences}, 54(2), 2001, pp.~143--164.

\bibitem[Kronecker 1901]{Kronecker1901}
L. Kronecker,
\textit{Vorlesungen \"uber Zahlentheorie}
Vol. 1, Teubner, Leipzig, 1901.

\bibitem[Kronecker 1903]{Kronecker1903}
L. Kronecker,
\textit{Vorlesungen \"uber die Theorie der Determinanten}
Teubner, Leipzig, 1903.

\bibitem[Lagrange 1795]{Lagrange1795}
J.L. Lagrange,
Le\c{c}on cinqui\`eme~: sur l'usage des courbes dans 
la solution des probl\`emes,
in \textit{S\'eances des \'Ecoles Normales recueillies 
par les st\'enographes et 
revues par les professeurs}, 
Reynier, Paris, 1795.

\bibitem[Lang 1970]{Lang1970}
S. Lang, \textit{Introduction to linear algebra},
Addison-Wesley, Reading, MA, 1970.

\bibitem[Lebesgue 1937]{Lebesgue1958}
H. Lebesgue
L'{\oe}uvre math\'ematique de Vandermonde,
in \textit{Notices d'Histoire des Math\'ematiques},
Universit\'e de Gen\`eve 1958, pp.~18--39.

\bibitem[Mansion 1880]{Mansion1880}
P. Mansion,
\textit{\'El\'ements de la th\'eorie des d\'eterminants},
Mons, 1880.

\bibitem[Marchand 1886]{Marchand1886}
E. Marchand, Sur le changement de variables,
\textit{Annales scientifiques de l'\'E.N.S.} 
3\textsuperscript{e} s\'erie, tome 3, 1886, pp.~137--138.

\bibitem[Merton 1968]{Merton1968}
R.K. Merton, The Matthew effect in science: the reward and
communication systems of science are considered,
\textit{Science}, 159(3810), 1968, pp.~56--63.

%\bibitem[Mikami 1914]{Mikami1914}
%Y. Mikami,
%On the Japanese theory of determinants,
%\textit{Isis}, 2(1), 1914, pp.~9--36.

\bibitem[Muir 1906]{Muir1906}
T. Muir,
\textit{The theory of determinants in the historical 
order of development, Volume I}
Mac Millan, London, 1906.

\bibitem[Muir 1911]{Muir1911}
T. Muir,
\textit{The theory of determinants in the historical 
order of development, Volume II}
Mac Millan, London, 1911.

\bibitem[Neuberg 1866]{Neuberg1866}
M. Neuberg, Question 590,
\textit{Nouvelles Annales de Math\'ema\-tiques},
2\textsuperscript{e} s\'erie, tome 5, 1866, p.~511--525.

\bibitem[Neumann 2007]{Neumann2007}
From Euler through Vandermonde to Gauss,
in R.E. Bradley \& C.E. Sandifer eds.
\textit{Leonhard Euler: life and legacy},
Elsevier, Amsterdam, 2007.

\bibitem[Newton 1687]{Newton1687}
I. Newton,
\textit{Philosophi{\ae} Naturalis Principia Mathematica},
London 1687.

\bibitem[Nielsen 1929]{Nielsen1929}
N. Nielsen,
\textit{G\'eom\`etres fran\c{c}ais sous la R\'evolution},
Levin \& Munksgaard, Copenhague, 1929.

\bibitem[Nuessel 2011]{Nuessel2011}
A note on the name of mathematical problems and
puzzles,
\textit{Names: J. Onomastics}, 50(1), 2011, pp.~57--64.

\bibitem[Pascal 1897]{Pascal1897}
E. Pascal,
\textit{Determinanti, teoria ed applicazioni},
Hoepli, Milano, 1897.

\bibitem[P\'olya \& Szeg\H{o} 1945]{Polya1945}
G. P\'olya and G. Szeg\H{o},
\textit{Problems and theorems in analysis},
Springer, New York, 1998.

\bibitem[Un Professeur 1860]{Anonymous1860}
Un Professeur,
Solution de la question 515,
\textit{Nouvelles Annales de Math\'ematiques}, 
1(19) 1860, pp.~181--183.

\bibitem[Przytycki 1992]{Przytycki1992}
J.H. Przytycki,
History of the knot theory from Vandermonde
to Jones,
in \textit{\textsc{xxiv}\textsuperscript{th} National Congress
of the Mexican Mathematical Society}, 
Mexico City, 1992, pp.~173--185.

\bibitem[Prouhet 1856]{Prouhet1856}
E. Prouhet, Notes sur quelques identit\'es,
\textit{Nouvelles Annales de Math\'ematiques},
1\textsuperscript{e} s\'erie, tome 15, 1856, pp.~86--91.

\bibitem[Remeslennikov 1993]{Remeslennikov1993}
V.N. Remeslennikov,
Vandermonde determinant
in M. Hazewinkel ed. \textit{Encyclopedia of Mathematics} vol. 9,
Kluwer, Dordrecht, 1993, p.~363.

\bibitem[Salmon 1859]{Salmon1859}
G. Salmon,
\textit{Lessons introductory to the modern higher algebra},
Hodges, Smith \& Co., Dublin, 1859.

\bibitem[Santos 2011]{Santos2011}
D.A. Santos,
\textit{Probability: an introduction},
Jones \& Bartlett, Sudbury MA, 2011.

\bibitem[Scott 1880]{Scott1880}
R.F. Scott,
\textit{A treatise on the theory of determinants,
and their applications to geometry and analysis},
Cambridge, 1880.

\bibitem[Scharf \emph{et al.} 1994]{Scharf1994}
T. Scharf, J.Y. Thibon, and B.G. Wybourne,
Powers of the Vandermonde determinant and the quantum Hall effect,
\textit{J. Phys. A: Math. Gen}, 27, 1994, pp.~4211--4219.

\bibitem[Small 2004]{Small2004}
H. Small,
On the shoulders of {Robert Merton}: towards a normative theory of citation,
\textit{Scientometrics}, 60(1), 2004, pp.~71--79.

\bibitem[Smith 1980]{Smith1980}
J.D.H. Smith,
In defense of eponymy,
\textit{Math Intelligencer}, 3(2), 1980, pp.~89--90.

\bibitem[Spottiswoode 1851]{Spottiswoode1851}
W. Spottiswoode,
\textit{Elementary theorems relating to determinants},
London, 1851.

\bibitem[Stedall 2011]{Stedall2011}
J. Stedall,
\textit{From Cardano's great art to Lagrange's reflections:
filling a gap in the history of algebra},
Heritage of European Mathematics, European Mathematical Society, 2011.

\bibitem[Stigler 1999]{Stigler1999}
S.M. Stigler,
\textit{Statistics on the table: the history of statistical concepts
  and methods},
Harvard University Press, Cambridge MA, 1999.

\bibitem[Suarez \& Gasc\'o 1882]{Suarez1882}
A. Suarez \& L. Gasc\'o,
\textit{Lecciones de combinatoria con las determinantes 
y sus principales aplicaciones}
Alufre, Valencia, 1882.

\bibitem[Sullivan 1997]{Sullivan1997}
C.R. Sullivan, The first chair of political economy in France:
Alexandre Vandermonde and the \emph{Principles} of Sir James Steuart
at the Ecole normale of the year \textsc{iii},
\textit{Franch Historical Studies}, 20(4), 1997, pp.~635--664.

\bibitem[Sylvester 1840]{Sylvester1840} 
J.J. Sylvester,
On derivation of coexistence: Part I. Being the theory of
  simultaneous simple homogeneous equations,
\textit{Philosophical Magazine}, vol.~\textsc{xvi}, 1840, pp.~37--43. 

\bibitem[Tee 1993]{Tee1993}
G.T. Tee,
Integer sums of recurring series,
\textit{New Zealand J. of Math.}, 22, 1993, pp.~85--100.

\bibitem[Terquem 1846]{Terquem1846}
O. Terquem,
Notice sur l'\'elimination,
\textit{Nouvelles Annales de Math\'ematiques}, 1\textsuperscript{e} 
s\'erie, tome 5, 1846, pp.~153--162.

\bibitem[Trudi 1862]{Trudi1862}
N. Trudi,
\textit{Teoria de' determinanti e loro applicazioni},
Pellerano, Napoli, 1862.

\bibitem[Utts \& Heckard 2004]{Utts2004}
J.M. Utts \& R.F. Heckard,
\textit{Mind on Statistics},
Thomson, Belmont, CA, 2004.

\bibitem[Vandermonde 1770]{Vandermonde1770}
A.T. Vandermonde,
M\'emoire sur la r\'esolution des
\'equa\-tions,
in 
\textit{Histoire de l'Acad\'emie royale des sciences
avec les m\'emoires de
math\'ematiques et de physique pour la 
m\^eme ann\'ee tir\'es des
registres de cette acad\'emie. 
Ann\'ee \textsc{mdcclxxi}}, Paris 1774,
pp.~365--413.

\bibitem[Vandermonde 1771]{Vandermonde1771}
A.T. Vandermonde,
M\'emoire sur l'\'elimination,
in 
\textit{Histoire de l'Acad\'emie royale des sciences
avec les m\'emoires de
math\'ematiques et de physique pour la 
m\^eme ann\'ee tir\'es des
registres de cette acad\'emie. 
Ann\'ee \textsc{mdcclxxi} seconde partie}, 
Paris 1776,
pp.~516--532.

\bibitem[Vandermonde 1772a]{Vandermonde1772a}
A.T. Vandermonde,
M\'emoire sur les irrationnelles des diff\'erents ordres avec une
application au cercle,
in 
\textit{Histoire de l'Acad\'emie royale des sciences
avec les m\'emoires de
math\'ematiques et de physique pour la 
m\^eme ann\'ee tir\'es des
registres de cette acad\'emie. 
Ann\'ee \textsc{mdcclxxii} premi\`ere partie}, 
Paris 1775,
pp.~489--498.

\bibitem[Vandermonde 1772b]{Vandermonde1772b}
A.T. Vandermonde,
Remarques sur les probl\`emes de situation,
in 
\textit{Histoire de l'Acad\'emie royale des sciences
avec les m\'emoires de
math\'ematiques et de physique pour la 
m\^eme ann\'ee tir\'es des
registres de cette acad\'emie. 
Ann\'ee \textsc{mdcclxxii} seconde partie}, 
Paris 1776,
pp.~566--574.

\bibitem[van der Waerden 1985]{Waerden1985}
B.L. van der Waerden,
\textit{A history of algebra: from al
Khaw\={a}rizm\={\i} to Emmy Noether},
Springer, New York 1985.

\bibitem[Vein \& Dale 1999]{Vein1999}
R. Vein \& P. Dale,
\textit{Determinants and their applications
in mathematical physics},
Applied Mathematical Sciences 134, Springer,
New York, 1999.

\bibitem[Wedderburn 1934]{Wedderburn1934}
J.H. Wedderburn,
\textit{Lectures on matrices},
AMS Colloquium publications, vol. 17, 1934.

\bibitem[Weill 1888]{Weill1888}
G. Weill,
Sur une forme du d\'eterminant de Vandermonde;
\textit{Nouvelles Annales de Math\'ematiques}, (3) vii, pp.~427--429, 1888.

\bibitem[Weld 1906]{Weld1906}
L.G. Weld,
\textit{Determinants},
Wiley, New York, 1906.

\bibitem[Whittacker \& Robinson 1924]{Whittacker1924}
E.T. Whittacker \& G. Robinson,
\textit{The calculus of observations},
Old Bailey, London 1924.
\end{thebibliography}

%
\end{document}